\documentclass[12pt,oneside,english]{amsart}
\usepackage[latin1]{inputenc}
\usepackage{amssymb}

\makeatletter
\newtheorem{theorem}{Theorem}
\newtheorem{lemma}{Lemma}

\newtheorem{proposition}{Proposition}
\newtheorem{remark}{Remark}

\numberwithin{equation}{section}
\usepackage{babel}

\makeatother
\begin{document}
\baselineskip=17pt
\title[On a sequence of polynomials  ]{On a sequence of polynomials with hypothetically integer coefficients}

\author{Vladimir Shevelev and Peter J. C. Moses}
\address{e-mail: shevelev@bgu.ac.il,\newline \indent \indent \indent  mows@mopar.freeserve.co.uk}
\subjclass{26C15}

\begin{abstract}
 The first author introduced a sequence of polynomials (\cite{8}, sequence A174531) defined recursively.  One of the main results of this study is proof of the integrality of its coefficients.
\end{abstract}

\maketitle
\section{Introduction}
In point of fact, there are only a few examples of sequences known where the question of the integrality of the terms is a difficult problem. In 1989, Somos \cite{9} posed a problem on the integrality of sequences depending on parameter $k\geq4$ which are defined by the recursion
\begin{equation}\label{1.1}
a_n=\frac{\sum_{j=1}^{\lfloor \frac{k}{2}\rfloor} a_{n-j}a_{n-(k-j)}} {a_{n-k} }, \enskip n\geq k\geq4,
\end{equation}
with the initial conditions $a_i=1, \enskip i=1,...,k-1.$\newline
 Gale \cite{3} proved the integrality of Somos sequences when $k=4$ and 5, attributing a proof to Malouf \cite{4}. Hickerson and Stanley (see \cite{6}) independently proved the integrality of the $k=6$ case in unpublished work and Fomin and Zelevinsky (2002) gave the first published proof. Finally, Lotto (1990) gave an unpublished proof for the $k=7$ case. These are sequences  A006720-A006723 in \cite{8}. It is interesting that, for $k\geq8,$ the property of integrality disappears (see sequence A030127 in \cite{8}). In connection with this, note that in the so-called G\"{o}bel's sequence (\cite{11}) defined by the recursion
\begin{equation}\label{1.2}
x_n=\frac{1}{n}(1+\sum_{i=0}^{n-1}x_i^2), \enskip n\geq1,\enskip x_0=1,
\end{equation}
the first non-integer term is $x_{43}=5.4093\times 10^{178485291567}.$ \newline
\indent In this paper we study the Shevelev sequence of polynomials $\{P_n(x)\}_{n\geq 1}$ that are defined by the following recursion $P_1=1,\enskip P_2=1,$ and, for $n\geq2,$

$$4(2x+n)P_{n+1}(x) =2(x+n)P_n(x)+$$
\begin{equation}\label{1.3}
(2x+n)P_n(x+1) +(4x+n)l_n(x), \enskip if \enskip n \enskip is \enskip odd,
\end{equation}\newline
 $$4P_{n+1}(x) = 4(x+n)P_n(x)+$$
\begin{equation}\label{1.4}
 2(2x+n+1)P_n(x+1)+(4x+n)l_{n-1}(x), \enskip if \enskip n \enskip is \enskip even,
\end{equation}
where
\begin{equation}\label{1.5}
l_n(x)=(x+\frac{n-1}{2})(x+\frac{n-3}{2})\cdot...\cdot(x+1).
\end{equation}

The first few polynomials are the following (\cite{8}, sequence A174531):
$$P_1=1,$$
$$P_2=1,$$
$$P_3=3x+4,$$
$$P_4=2x+4,$$
$$P_5=5x^2+25x+32,$$
$$P_6=3x^2+19x+32,$$
$$P_7=7x^3+77x^2+294x+384,$$
$$P_8=4x^3+52x^2+240x+384,$$
$$P_9=9x^4+174x^3+1323x^2+4614x+6144,$$
$$P_{10}=5x^4+110x^3+\enskip967x^2+3934x+6144,$$
$$P_{11}=11x^5+330x^4+4169x^3+27258x^2+90992x+122880,$$
$$P_{12}=\enskip 6x^5+200x^4+2842x^3+21040x^2+79832x+122880.$$
According to our observations, the following conjectures are natural.\newline
1) The coefficients of all the polynomials are integers. Moreover, the greatest common divisor of all coefficients is ${n}/{rad(n)},$ where $rad(n)=\prod_{p|n}p;$ \newline
2) $P_n(0)=4^{\lfloor\frac{n-1}{2}\rfloor}\lfloor \frac{n-1}{2}\rfloor!;$ \newline
3) For even $n$, $P_n(1)=(2^n-1)(\frac {n}{2})!/(n+1),$ and for odd $n$,  $P_n(1)=(2^n-1)(\frac{n-1}{2})!;$\newline
4) $P_n(x)$ has a real rational root if and only if either $n=3$ or $n\equiv0\pmod4.$ In the latter case, such a unique root is $-\frac{n}{2};$ \newline
5) Coefficients of $x^k$ increase when $k$ decreases;\newline
6) If $n$ is even, then the coefficients of $P_n$ do not exceed the corresponding coefficients of $P_{n-1}$ and the equality holds only for the last ones; moreover, the ratios of coefficients of $x^k$ of polynomials $P_{n-1}$ and $P_n$ monotonically decrease to 1 when $k$ decreases; \newline
7) All coefficients of $P_n,$ except of the last one, are multiple of $n$ if and only if $n$ is prime.\newline
\newpage
The main results of our paper consist of the following two theorems.
\begin{theorem}\label{t1} \upshape(Explicit formula for $P_n(k))$\slshape \enskip For integer $x=k,$ we have
\begin{equation}\label{1.6}
P_n(k) =
\end{equation}
$$\begin{cases}
(\binom{(n-1)/2+k-1}{k-1}/\binom{n+2k-2}{k-1})(\frac {n-1}{2})!\enskip T_n(k),\;\;if\;\;n\geq1\;\;is\;\;odd,\\\\ (\binom{n/2+k-1}{k} / \binom{n+2k-1}{k})(\frac{n}{2}-1)!\enskip T_n(k),\;\;if\;\; n\geq2\;\; is\;\; even,\end{cases}$$

\begin{equation}\label{1.7}
=2^{-(\lfloor\frac{n}{2}\rfloor+k-1)}\frac{(n+k-1)!}{(2\lfloor\frac{n}{2}\rfloor+2k-1)!!}T_n(k),
\end{equation}
where
\begin{equation}\label{1.8}
T_n(k)=\sum_{i=1}^n 2^{i-1}\binom{n+2k-i-1}{k-1}.
\end{equation}
\end{theorem}
Using Theorem \ref{t1}, we prove Conjectures 2)-3) and the following main result.

\begin{theorem}\label{t2} For $n\geq1,$
 $P_n(x)$ is a polynomial of degree $\lfloor\frac{n-1}{2}\rfloor$ with integer coefficients.
\end{theorem}
Nevertheless, the subtle second part of conjecture 1) remains open.
\section{Representation of $P_n(k)$ via a polynomial in $n$ of degree $k-1$ with integer coefficients}
\begin{theorem}\label{t4} For integer $k\geq1,\enskip n\geq1,$ the following recursion holds
\begin{equation}\label{2.1}
P_n(k)=c_n(k)(2^{n+k-1}-\frac {R_{k}(n)}{(2k-2)!!}),
\end{equation}
where $R_k(n)$ is a polynomial in $n$ of degree $k-1$ with integer coefficients and
\begin{equation}\label{2.2}
c_n(k)=\begin{cases}
(\frac {n-1}{2})!\prod_{i=1}^{k-1}\frac {n+i}{n+2i},\;\;if\;\;n\;\;is\;\;odd,\\\\ \frac{1}{2}(\frac{n}{2}-1)!\prod_{i=0}^{k-1}\frac {n+i}{n+2i+1},\;\;if\;\; n\;\; is\;\; even,\end{cases}
\end{equation}
\end{theorem}
 \bfseries Proof. \mdseries
 Write (\ref{1.3})-(\ref{1.4}) in the form
 $$P_n(k+1)=-\frac{2f}{g}P_n(k)+4P_{n+1}(k)-$$
\begin{equation}\label{2.3}
\frac{h}{g}(\frac {n-1}{2})!\binom{\frac{g-1}{2}}{k},\enskip if \enskip n\equiv1\pmod2;
\end{equation}
 $$P_n(k+1)=-\frac{2f}{g+1}P_n(k)+\frac{2}{g+1}P_{n+1}(k)-$$
 \newpage
\begin{equation}\label{2.4}
\frac{h}{2(g+1)}(\frac {n}{2}-1)!\binom{\frac{g}{2}-1}{k},\enskip if \enskip n\equiv0\pmod2,
\end{equation}
where $f=n+k, \enskip g=n+2k, \enskip h=n+4k.$\newline
\indent Let $n$ be odd. We use induction over $k$. For $k=1,$ (\ref{2.1}) gives
\begin{equation}\label{2.5}
R_1(n)=2^{n}-\frac{P_n(1)}{c_n(1)}=Const(k).
\end{equation}
Thus the base of induction is valid. Suppose the theorem is true for some value of $k.$ Then, using this supposition and (\ref{2.1})-(\ref{2.4}), we have
$$ P_n(k+1)=$$$$-\frac{2f}{g}(\frac {n-1}{2})!(2^{n+k-1}-\frac {R_{k}(n)}{(2k-2)!!})\prod_{i=1}^{k-1}\frac {n+i}{n+2i}+$$ $$2(\frac {n-1}{2})!(2^{n+k}-\frac {R_{k}(n+1)}{(2k-2)!!})\prod_{i=0}^{k-1}\frac {n+i+1}{n+2i+2}-$$
$$\frac{h}{g}(\frac {n-1}{2})!\frac{\frac{g-1}{2}\frac{g-3}{2}\cdot...\cdot\frac{n+1}{2}}{k!}.  $$
Note that
$$ \frac {f}{g}\prod_{i=0}^{k-1}\frac {n+i}{n+2i}=\prod_{j=1}^{k}\frac {n+j}{n+2j}=\prod_{i=0}^{k-1}\frac {n+i+1}{n+2i+2}.$$
Therefore,
$$ P_n(k+1)=(\frac {n-1}{2})!(-2^{n+k}+\frac {2R_{k}(n)}{(2k-2)!!}+2^{n+k+1}-$$$$\frac {2R_{k}(n+1)}{(2k-2)!!}-\frac{h}{g}\frac{\frac{g-1}{2}\frac{g-3}{2}\cdot...\cdot\frac{n+1}{2}}{k!}\prod_{j=1}^{k}\frac {n+2j}{n+j})\prod_{j=1}^{k}\frac {n+j}{n+2j}.$$
Here we note that
$$(g-1)(g-3)\cdot...\cdot(n+1)\prod_{j=1}^{k}\frac {n+2j}{n+j}=(n+2k)_k,$$
where $(x)_k$ is a falling factorial.
Hence
$$ P_n(k+1)=c_n(k+1)(2^{n+k}-2\frac {R_{k}(n+1)-R_k(n)}{(2k-2)!!}-$$$$ \frac{4k+n}{(2k)!!}(n+2k-1)_{k-1})=c_n(k+1)(2^{n+k}-\frac {R_{k+1}(n)}{(2k)!!}),$$
where
\begin{equation}\label{2.6}
R_{k+1}(n)= 4k(R_{k}(n+1)-R_k(n))+(4k+n)(n+2k-1)_{k-1}.
\end{equation}
Since, by the inductive supposition, $R_k(n)$ is a polynomial of degree $k-1$ with integer coefficients, then, by (\ref{2.6}), $R_{k+1}(n)$ is a polynomial of degree $k$ with integer coefficients.\newpage
\indent Note that the case of even $n$ is considered quite analogously, obtaining the \slshape same\upshape\enskip formula (\ref{2.6}). $ \blacksquare$\newline
\indent Put in (\ref{2.1})-(\ref{2.2}) $n=1$. Then, for $k\geq1$ we have
$$(2^k-\frac{R_k(1)}{(2k-2)!!})\frac{k!}{(2k-1)!!}=1,$$
whence
\begin{equation}\label{2.7}
R_{k}(1)=(k-1)!(2^{2k-1}-\binom{2k-1}{k}).
\end{equation}
In particular, $R_1(1)=1$ and, since $R_1(n)$ is of degree 0, then $R_1(n)=1.$
Further, we find polynomials $R_k(n)$ using the recursion (\ref{2.6}). The first polynomials $R_k(n)$ are
$$R_1(n)=1,$$
$$R_2(n)=n+4,$$
$$R_3(n)=n^2+11n+32,$$
$$R_4(n)=n^3+21n^2+152n+384,$$
$$R_5(n)=n^4+34n^3+443n^2+2642n+6144,$$
$$R_6(n)=n^5+50n^4+1015n^3+10510n^2+55864n+122880. $$
\section{Proof of Conjectures 2) and 3)}
We start with proof of Conjecture 3) for $P_n(1).$  Note that, since $R_1(n)=1,$  then from (\ref{2.5}) we find
\begin{equation}\label{3.1}
P_n(1)=c_n(1)(2^n-1).
\end{equation}
Besides, by (\ref{2.2}), we have
\begin{equation}\label{3.2}
c_n(1)=\begin{cases}
(\frac {n-1}{2})!,\;\;if\;\;n\;\;is\;\;odd,\\\\ \frac{1}{2}(\frac{n}{2}-1)!\frac {n}{n+1}=(\frac{n}{2})!/(n+1),\;\;if\;\; n\;\; is\;\; even,\end{cases}
\end{equation}
and Conjecture 3) follows. \newline
Let us prove now Conjecture 2). Note that (\ref{2.3})-(\ref{2.4}), as (\ref{1.3})-(\ref{1.4}), is valid for every nonnegative $k.$ For $k=0$ and odd $n\geq1,$ (\ref{2.3}) gives
$$P_n(1)=-2P_n(0)+4P_{n+1}(0)-(\frac{n-1}{2})!,$$
or, using (\ref{3.1})-(\ref{3.2}), we have
$$4P_{n+1}(0)-2P_n(0)=2^n(\frac{n-1}{2})!$$
Analogously, for $k=0$ and even $n\geq1,$ from (\ref{2.4}) and (\ref{3.1})-(\ref{3.2}) we find

$$P_{n+1}(0)-nP_n(0)=2^{n-1}(\frac{n}{2})!$$
\newpage
Thus
\begin{equation}\label{3.3}
P_{n+1}(0)=\begin{cases}
\frac {1}{2}P_n(0)+2^{n-2}(\frac{n-1}{2})!,\;\;if\;\;n\;\;is\;\;odd,\\\\ nP_n(0)+2^{n-1}(\frac{n}{2})!,\;\;if\;\; n\;\; is\;\; even\end{cases}
\end{equation}
with $P_1(0)=1,\enskip P_2(0)=1.$ Since the difference equation
$$y(n+1)=\begin{cases}
\frac {1}{2}y(n)+2^{n-2}(\frac{n-1}{2})!,\;\;if\;\;n\;\;is\;\;odd,\\\\ ny(n)+2^{n-1}(\frac{n}{2})!,\;\;if\;\; n\;\; is\;\; even\end{cases}$$
with the initials $y(1)=1,\enskip y(2)=1$ has an unique solution, then it is sufficient to verify that $y(n)=P_n(0)=4^{\lfloor\frac{n-1}{2}\rfloor}\lfloor \frac{n-1}{2}\rfloor!$ is a solution. $\blacksquare$\newline
\section{Explicit formula for $R_k(n)$}
Since from (\ref{2.6})
$$4kR_{k}(n+1)=4kR_k(n)+$$
\begin{equation}\label{4.1}
R_{k+1}(n)-(4k+n)(n+2k-1)_{k-1},
\end{equation}
we have a recursion in $n$ for $R_k(n)$ given by (\ref{2.7}) and (\ref{4.1}).\newline
Our aim in this section is to find a generalization of (\ref{2.7}) for an arbitrary integer $n\geq1.$
Note that we can write (\ref{2.7}) in the form
\begin{equation}\label{4.2}
R_k(1)=2(k-1)!4^{k-1}-\frac{(2k-1)!}{k!}.
\end{equation}
Using (\ref{4.1}) and (\ref{2.7}), after some transformations, we find
\begin{equation}\label{4.3}
R_k(2)=2^2(k-1)!4^{k-1}-2\frac{(2k-1)!}{k!}-\frac{(2k)!}{(k+1)!}.
\end{equation}
The regularity is fixed in the following theorem.
\begin{theorem}\label{t5} For integer $k\geq1,\enskip n\geq1,$ we have
\begin{equation}\label{4.4}
R_k(n)=2^n(k-1)!4^{k-1}-\sum_{i=1}^n 2^{n-i}\frac {(2k+i-2)!}{(k+i-1)!}.
\end{equation}
\end{theorem}
 \bfseries Proof. \mdseries Taking into account that $\frac {(2k+i-2)!}{(k+i-1)!}=\binom{2k+i-2}{k-1}(k-1)!,$ we prove (\ref{4.4}) in the following equivalent form:
\begin{equation}\label{4.5}
R_k(n)=2^n(k-1)!(4^{k-1}-\sum_{i=1}^n 2^{-i}\binom{2k+i-2}{k-1}).
\end{equation}
We use induction over $n.$ Suppose that (\ref{4.5}) is valid for a some value of $n$ and an arbitrary integer $k\geq1.$ Then, by (\ref{4.1}), we have
$$R_k(n+1)=2^n(k-1)!(4^{k-1}-\sum_{i=1}^n 2^{-i}\binom{2k+i-2}{k-1})+$$
\newpage
$$2^{n-2}(k-1)!(4^{k}-\sum_{i=1}^n 2^{-i}\binom{2k+i}{k})-\frac{4k+n}{4k}(n+2k-1)_{k-1}=$$
$$2^n(k-1)!(4^{k-1}-\sum_{i=1}^n 2^{-i}\binom{2k+i-2}{k-1})+$$
$$2^n(k-1)!(4^{k-1}-\sum_{i=1}^n 2^{-i-2}\binom{2k+i}{k})-$$
$$\frac{(n+2k-1)!}{(n+k)!}-\frac{n}{4k}\frac{(n+2k-1)!}{(n+k)!}.$$
Thus we should prove the identity
$$2^{n+1}(k-1)!4^{k-1}-2^n(k-1)!\sum_{i=1}^n 2^{-i}\binom{2k+i-2}{k-1}-$$
$$2^{n-2}(k-1)!\sum_{i=1}^n 2^{-i}\binom{2k+i}{k}-\frac{n+4k}{4k}\frac{(n+2k-1)!}{(n+k)!}=  $$
$$2^{n+1}(k-1)!(4^{k-1}-\sum_{i=1}^{n+1} 2^{-i}\binom{2k+i-2}{k-1}),$$
which is easily reduced to the identity
$$4\sum_{i=1}^{n} 2^{-i}\binom{2k+i-2}{k-1}-\sum_{i=1}^{n} 2^{-i}\binom{2k+i}{k}=$$
$$2^{-n}\frac{n+4k}{4k}\frac{(n+2k-1)!}{(n+k)!}-4\cdot2^{-n}\binom{2k+n-1}{k-1}.$$
Note that, the right hand part is $\frac{n}{k2^n}\binom{2k+n-1}{k-1}.$ Therefore, it is left to prove the identity
\begin{equation}\label{4.6}
4\sum_{i=1}^{n} 2^{-i}\binom{2k+i-2}{k-1}-\sum_{i=1}^{n} 2^{-i}\binom{2k+i}{k}=\frac{n}{k2^n}\binom{2k+n-1}{k-1}.
\end{equation}
Since this is trivially satisfied for $n=0,$ then it is sufficient to verify the equality of the first differences of the left and the right hand parts, which is reduced to the identity
$$ 2(n+2k-1)\binom{2k+n-2}{k-1}= {n}\binom{2k+n-1}{k-1}+k\binom{2k+n}{k},$$
which is verified directly. $\blacksquare$\newpage
\section{Proof of Theorem \ref{t1}}
Now we are able to prove Theorem \ref{t1}. According to (\ref{1.7}), we have
$$T_n(k)=\sum_{i=1}^n 2^{i-1}\binom{n+2k-i-1}{k-1}=$$
\begin{equation}\label{5.1}
\sum_{j=1}^n 2^{n-j}\binom{2k+j-2}{k-1}.
\end{equation}
Hence, by (\ref{4.5}), we find
$$ R_k(n)=2^n(k-1)!(4^{k-1}-2^{-n}T_n(k))=$$
\begin{equation}\label{5.2}
(k-1)!(2^{n+2k-2}-T_n(k)).
\end{equation}
Now from (\ref{2.1}) and (\ref{5.2}) we have
\begin{equation}\label{5.3}
P_n(k)=2^{-(k-1)}c_n(k)T_n(k).
\end{equation}
Let $n$ be odd.
Note that, by (\ref{2.2}),
$$2^{-(k-1)}c_n(k)=$$
 $$2^{-(k-1)}(\frac{n-1}{2})!\frac {(n+k-1)(n+k-2)...(n+1)}{(n+2k-2)(n+2k-4)...(n+2)}=$$
\begin{equation}\label{5.4}
 2^{-(k-1)}(\frac{n-1}{2})!\frac{(n+k-1)!n!!}{n!(n+2k-2)!!}.
\end{equation}
 Taking into account that
\begin{equation}\label{5.5}
  n!!=\frac{n!}{(n-1)!!}=\frac{n!}{2^\frac{n-1}{2}(\frac{n-1}{2})!},
\end{equation}
we find from (\ref{5.4})
$$ 2^{-(k-1)}c_n(k)=\frac {(n+k-1)!(\frac{n-1}{2}+k-1)!}{(n+2k-2)!}=$$
$$\frac {(\frac{n-1}{2}+k-1)!} {(k-1)!\binom{n+2k-2}{k-1}}=\frac {\binom{\frac{n-1}{2}+k-1}{k-1}}{\binom{n+2k-2}{k-1}}(\frac{n-1}{2})!$$
and ({1.6}) follows from (\ref{5.3}). Furthermore, since by (\ref{5.5}) $\frac{n!!(\frac{n-1}{2})!}{n!}=2^{-\frac{n-1}{2}},$ then from (\ref{5.3})-(\ref{5.4}) we find
$$P_n(k)=2^{-(\frac{n-1}{2}+k-1)}\frac{(n+k-1)!}{(n+2k-2)!!}T_n(k)$$
that corresponds to ({1.7}) in the case of odd $n.$ The case of even $n$ is considered quite analogously. $\blacksquare$\newpage
\section{Bisection of sequence $\{P_n(x)\}$}
Note that $T_n(k)$ (\ref{1.8}) has rather a simple structure, which allows us to find different relations for it. Using (\ref{1.6}), this, in turn, allows us to find recursion relations for $P_n(x)$ which are simpler than the basis recursion (\ref{1.3})-(\ref{1.4}). We start with the following simple recursions for $T_n(k).$
\begin{lemma}\label{L1}
\begin{equation}\label{6.1}
T_n(k)-2T_{n-1}(k)=\binom{n+2k-2}{k-1},\enskip k\geq1;
\end{equation}
\begin{equation}\label{6.2}
T_n(k)-4T_{n-2}(k)=\binom{n+2k-2}{k-1}+2\binom{n+2k-3}{k-1},\enskip k\geq2.
\end{equation}
\end{lemma}
\bfseries Proof. \mdseries By (\ref{1.8}), we have
$$T_n(k)-2T_{n-1}(k)=$$ $$\sum_{i=1}^n 2^{i-1}\binom{n+2k-i-1}{k-1}-\sum_{j=1}^{n-1} 2^{j}\binom{n+2k-j-2}{k-1}=$$
$$ \sum_{i=1}^n 2^{i-1}\binom{n+2k-i-1}{k-1}-\sum_{i=2}^{n} 2^{i-1}\binom{n+2k-i-1}{k-1}$$
and (\ref{6.1}) follows; (\ref{6.2}) is a simple corollary of (\ref{6.1}). $\blacksquare$
\begin{theorem}\label{t5} (Bisection)
If $n\geq3$ is odd, then
$$(2x+n-2)P_n(x)=2(x+n-1)(x+n-2)P_{n-2}(x)+$$
\begin{equation}\label{6.3}
(4x+3n-4)(x+\frac{n-1}{2}-1)(x+\frac{n-1}{2}-2)\cdot...\cdot x;
\end{equation}
if $n\geq4$ is even, then
$$(2x+n-1)P_n(x)=2(x+n-1)(x+n-2)P_{n-2}(x)+$$
\begin{equation}\label{6.4}
\frac{1}{2}(4x+3n-4)(x+\frac{n-2}{2}-1)(x+\frac{n-2}{2}-2)\cdot...\cdot x.
\end{equation}
\end{theorem}
 \bfseries Proof. \mdseries According to (\ref{1.6}), we have
 \begin{equation}\label{6.5}
T_n(k) =\begin{cases}
\binom{n+2k-2}{k-1}/(\binom{(n-1)/2+k-1}{k-1}(\frac {n-1}{2})!)\enskip P_n(k),\;\;if\;\;n\;\;is\;\;odd,\\\\ \binom{n+2k-1}{k}/(\binom{n/2+k-1}{k}(\frac{n}{2}-1)!)\enskip P_n(k),\;\;if\;\; n\;\; is\;\; even.\end{cases}
\end{equation}
Substituting this to (\ref{6.2}), after simple transformations, we obtain (\ref{6.3})-(\ref{6.4}), where $k$ is replaced by arbitrary $x.$ $\blacksquare$\newline
\indent Note that from ({6.3})-({6.4}), using a simple induction, we conclude that, for even $n\geq4,$ $P_n(x)$ is a polynomial of degree $\frac{n-2}{2},$ while, for odd $n\geq3,$  $P_n(x)$ is a polynomial of degree $\frac{n-1}{2}.$ However, a structure of formulas ({6.3})-({6.4}) does not allow us to prove that all coefficients of $P_n(x)$ are integer.
\newpage This will be done in the following section by the discovery of the special relationships with the required structure.

\section{Proof of Theorem 2}
\begin{lemma}\label{L2} For $n\geq1,$ we have
\begin{equation}\label{7.1}
T_n(k)-T_{n-2}(k+1)=\binom{n+2k-1}{k}.
\end{equation}
\end{lemma}
\bfseries Proof. \mdseries By (\ref{5.1}), we should prove that
$$\binom{2k+n-1}{k}=T_n(k)-T_{n-2}(k+1)=$$ $$\sum_{j=1}^n 2^{n-j}\binom{2k+j-2}{k-1}-\sum_{j=1}^{n-2} 2^{n-j-2}\binom{2k+j}{k}=$$
$$\sum_{j=1}^n 2^{n-j}\binom{2k+j-2}{k-1}-\sum_{i=1}^{n} 2^{n-i}\binom{2k+i-2}{k}+$$
$$2^{n-1}\binom{2k-1}{k}+2^{n-2}\binom{2k}{k}, $$
or
$$\sum_{j=1}^n 2^{-j}(\binom{2k+j-2}{k-1}-\binom{2k+j-2}{k})=$$
\begin{equation}\label{7.2}
2^{-n}\binom{2k+n-1}{k}-\frac{1}{2}\binom{2k-1}{k}-\frac{1}{4}\binom{2k}{k}.
\end{equation}
It is verified directly that ({7.2}) is valid for $n=1.$ Therefore, it is sufficient to verify that the first differences over $n$ of the left hand side and the right hand side coincide. The corresponding identity
$$2^{-n}(\binom{2k+n-2}{k-1}-\binom{2k+n-2}{k})=$$ $$2^{-n}\binom{2k+n-1}{k}-2^{-n+1}\binom{2k+n-2}{k}$$
reduces to the equality $\binom{2k+n-2}{k-1}+\binom{2k+n-2}{k}=\binom{2k+n-1}{k}.\blacksquare$ \newline
\indent Now we are able to complete proof of Theorem 2.
Considering even $n\geq4,$ by ({6.5}), we obtain the following relation for $P_n(k)$ corresponding to ({7.1}):
$$P_n(x)=(n+x-1)P_{n-2}(x+1)+$$
\begin{equation}\label{7.3}
(x+\frac{n}{2}-1)(x+\frac{n}{2}-2)...(x+1).
\end{equation}
On the other hand, using (\ref{6.1}), for odd $n\geq3,$ we obtain the following relation
$$P_n(x)=2(x+n-1)P_{n-1}(x)+$$
\newpage
\begin{equation}\label{7.4}
(x+\frac{n-1}{2}-1)(x+\frac{n-1}{2}-2)\cdot...\cdot x.
\end{equation}

\indent From ({7.3}), by a simple induction, we see that, for even $n\geq4,$ $P_n(x)$ is a polynomial with integer coefficients. Then from ({7.4}) we find that $P_n(x),$ for odd $n,$ is a polynomial with integer coefficients as well. $\blacksquare$

\section{Other relations}
Together with ({6.3})-({6.4}), ({7.3})-({7.4}) there exist many other relations for $P_n(x).$ All of them are corollaries of the corresponding relations for $T_n(k).$ Below we give a few pairs of some such relations.\newline
\indent As we saw, for odd $n\geq3,$ ({7.4}) follows from ({6.1}). Let us consider even $n\geq4.$ Then we obtain the second component of the following recursion
\begin{equation}\label{8.1}
P_n(x)=\begin{cases}2(x+n-1)P_{n-1}(x)+\\\\((x+n-1)P_{n-1}(x)+\end{cases}
\end{equation}

$$\begin{cases}
(x+\frac{n-1}{2}-1)(x+\frac{n-1}{2}-2)\cdot...\cdot x,\;\;if\;\;n\geq3\;\;is\;\;odd,\\\\(x+\frac{n}{2}-1)(x+\frac{n}{2}-2)\cdot...\cdot x)/(2x+n-1) ,\;\;if\;\; n\geq4\;\; is\;\; even.\end{cases}$$

\begin{lemma}\label{L3} For $n\geq1,\enskip k\geq1,$ we have
\begin{equation}\label{8.2}
T_n(k+1)=4T_{n}(k)-\frac{n}{k}\binom{n+2k-1}{k-1}.
\end{equation}
\end{lemma}
\bfseries Proof. \mdseries By ({7.1}),({6.2}), we have
$$T_n(k+1)=T_{n+2}(k)-\binom{n+2k+1}{k}=$$ $$4T_n(k)+\binom {n+2k}{k-1}+2\binom {n+2k-1}{k-1}-\binom{n+2k+1}{k}.$$
It is left to note that
$$\binom {n+2k}{k-1}+2\binom {n+2k-1}{k-1}-\binom{n+2k+1}{k}=-\frac{n}{k}\binom{n+2k-1}{k-1}.$$
$\blacksquare$\newline
From Lemma \ref{L3} and ({6.5}) we find the following recursion
\begin{equation}\label{8.3}
\begin{cases}(2x+n)P_{n}(x+1)=2(x+n)P_n(x)-\\\\(2x+n+1)P_{n}(x+1)=2(x+n)P_n(x)-\end{cases}
\end{equation}

$$\begin{cases}
{n}(x+\frac{n-1}{2})(x+\frac{n-1}{2}-1)\cdot...\cdot (x+1),\;\;if\;\;n\geq3\;\;is\;\;odd,\\\\ \frac{n}{2}(x+\frac{n}{2}-1)(x+\frac{n}{2}-2)\cdot...\cdot (x+1) ,\;\;if\;\; n\geq4\;\; is\;\; even.\end{cases}$$
\newpage
\begin{lemma}\label{L4} For $n\geq2,\enskip k\geq1,$ we have
\begin{equation}\label{8.4}
(n+k-1)(T_n(k)-4T_n(k-1))=n(T_{n-1}(k)-2T_n(k-1)).
\end{equation}

\end{lemma}
\bfseries Proof. \mdseries By ({8.2}),
\begin{equation}\label{8.5}
T_n(k)-4T_n(k-1)=-\frac{n}{k-1}\binom {n+2k-3}{k-2}.
\end{equation}
By ({6.1}),
$$T_n(k-1)=2T_{n-1}(k-1)+\binom {n+2k-4}{k-2}.$$
Therefore,
$$T_{n-1}(k)-2T_n(k-1)=T_{n-1}(k)-4T_{n-1}(k-1)-2\binom {n+2k-4}{k-2}.$$
Using again ({8.2}), we find
\begin{equation}\label{8.6}
T_{n-1}(k)-2T_n(k-1)=-(\frac{n-1}{k-1}+2)\binom {n+2k-4}{k-2}.
\end{equation}
Now the lemma follows from ({8.5})-({8.6}) since $(n+k-1)\binom {n+2k-3}{k-2}=(n+2k-3)\binom {n+2k-4}{k-2}. \blacksquare$\newline
\indent The passage from ({8.4}) to the corresponding formula for $P_n(x)$ in the case of odd $n\geq3$ unexpectedly leads to a very simple homogeneous relation
\begin{equation}\label{8.7}
 P_n(x)=P_n(x-1)+n P_{n-1}(x)
\end{equation}
which we use in Sections 9 and 12. The corresponding relation for even $n\geq4$ is
\begin{equation}\label{8.8}
 (2x+n-1)P_n(x)=(2x+n-2)P_n(x-1)+\frac{n}{2}P_{n-1}(x).
\end{equation}
\begin{lemma}\label{L5} For $n\geq1,\enskip k\geq2,$ we have
\begin{equation}\label{8.9}
2T_n(k)-T_{n-1}(k+1)=\binom{n+2k-1}{k}.
\end{equation}
\end{lemma}
\bfseries Proof. \mdseries By ({6.1}), we have
$$2T_n(k)-T_{n-1}(k+1)=$$ $$4T_{n-1}(k)+2\binom{n+2k-2}{k-1}-T_{n-1}(k+1).$$
Furthermore, by ({7.1}),
$$T_{n-1}(k+1)=T_{n+1}(k)-\binom{n+2k}{k}.$$
Hence,
$$2T_n(k)-T_{n-1}(k+1)=$$
\begin{equation}\label{8.10}
4T_{n-1}(k)-T_{n+1}(k)+2\binom{n+2k-2}{k-1}+\binom{n+2k}{k}.
\end{equation}
\newpage
Finally, by ({6.2}),
$$T_{n+1}(k)-4T_{n-1}(k)=\binom{n+2k-1}{k-1}+2\binom{n+2k-2}{k-1}$$
and the lemma follows from ({8.10}). $\blacksquare$\newline
Using Lemma \ref{L5} and ({6.5}), for even $n\geq4,$ we find
\begin{equation}\label{8.11}
2P_n(x)=P_{n-1}(x+1)+(x+\frac{n}{2}-1)(x+\frac{n}{2}-2)\cdot...\cdot(x+1),
\end{equation}
while, for odd $n\geq3,$
\begin{equation}\label{8.12}
P_n(x)=(2x+n)P_{n-1}(x+1)+(x+\frac{n-1}{2})(x+\frac{n-1}{2}-1)\cdot...\cdot(x+1).
\end{equation}

\begin{proposition}\label{P1}
For odd $n\geq3,$ we have
\begin{equation}\label{8.13}
P_n(k)\equiv P_n(0) \pmod n.
\end{equation}
\end{proposition}
\bfseries Proof. \mdseries From ({8.7}) we find
\begin{equation}\label{8.14}
\sum_{i=1}^k P_{n-1}(i)=(P_n(k)-P_n(0))/n,
\end{equation}
and the proposition follows. $\blacksquare$
\section{On coefficients of $P_n(x)$}

Using formulas (\ref{6.3})-(\ref{6.4}), we give a recursion for calculation of the coefficients of $P_n(x)$ with a fixed parity of $n.$ Let
\begin{equation}\label{9.1}
 P_n(x)=a_0(n){x}^{m}+a_1(n){x}^{m-1}+...+a_{m-1}(n){x}+a_m(n),
\end{equation}
where $m=\lfloor\frac{n-1}{2}\rfloor.$ We prove the following.
\begin{theorem}\label{t6} For $n\geq1,$ we have
\begin{equation}\label{9.2}
a_0(n)=\begin{cases}
n, \;\;if\;\;n\;\;is\;\;odd,\\\\\frac{n}{2},\;\;if\;\; n\;\; is\;\; even;\end{cases}
\end{equation}
$$a_1(n)=\begin{cases}
\frac{1}{24}(7n^3-12n^2+5n)\\\\\frac{1}{48}(7n^3-18n^2+8n)\end{cases}=$$
\begin{equation}\label{9.3}
\begin{cases}
\frac{1}{24}n(n-1)(7n-5), \;\;if\;\;n\;\;is\;\;odd,\\\\\frac{1}{48}n(n-2)(7n-4),\;\;if\;\; n\;\; is\;\; even.\end{cases}
\end{equation}
In general, for a fixed $i,$ $a_i(n)=U_i(n),$ if $n$ is odd, and $a_i(n)=V_i(n),$ if $n$ is even, where $U_i,\enskip V_i$ are polynomials in $n$ of degree $2i+1.$
\end{theorem}
\newpage
 \bfseries Proof. \mdseries
 1) Let $n$ be even. Then, using (\ref{6.4}), for integer $x$ and $m=\frac{n-2}{2},$ we have
 $$ (2x+n-1)(a_0(n)x^m+a_1(n)x^{m-1}+...)=$$ $$2(x+n-1)(x+n-2)(a_0(n-2)x^{m-1}+a_1(n-2)x^{m-2}+...)+$$
\begin{equation}\label{9.4}
 \frac{1}{2}(\frac{n-2}{2})!(4x+3n-4)\binom{x-1+\frac{n-2}{2}}{\frac{n-2}{2}}.
 \end{equation}
 Comparing the coefficient of $x^{m+1}$ in both hand sides, we find
 $$a_0(n)=a_0(n-2)+1,\enskip n\geq4,\enskip a_0(4)=2.$$
 Thus $a_0(6)=3, a_0(8)=4,..., a_0(n)=n/2.$  \newline
Furthermore, comparing the coefficient of $x^{m}$ in both hand sides in (\ref{9.4}), we have
$$2a_1(n)+(n-1)a_0(n)=2a_1(n-2)+2(2n-3)a_0(n-2)+$$
\begin{equation}\label{9.5}
Coef[x^m](\frac{1}{2}(4x+3n-4)(x+\frac{n-4}{2})(x+\frac{n-6}{2})\cdot...\cdot(x+1)x).
\end{equation}
Note that
$$Coef[x^m](\frac{1}{2}(4x+3n-4)(x+\frac{n-4}{2})(x+\frac{n-6}{2})\cdot...\cdot(x+1)x)=$$
$$\frac{3n-4}{2}+2(\frac{n-4}{2}+\frac{n-6}{2}+...+1)=$$
$$\frac{3n-4}{2}+\sum_{i=2}^{m}(n-2i)=\frac{n^2}{4}.$$
Therefore, by ({9.5}),
$$a_1(n)-a_1(n-2)=\frac{(2n-3)(n-2)}{2}-\frac{(n-1)n}{4}$$ $$+\frac{n^2}{8}=
\frac{7n^2-26n+24}{8}.$$
Hence
 $$a_1(n)=\sum_{i=4,6,...,n}(a_1(i)-a_1(i-2))=
  \frac{1}{8}\sum_{i=4,6,...,n}(7i^2-26i+24)=$$ $$\frac{1}{2}\sum_{j=2}^{n/2}(7j^2-13j+6)=\frac{1}{48}(7n^3-18n^2+8n).$$
 $\blacksquare$ \newline
 Finally, comparing the coefficient of $x^{m-i}$ in both hand sides of (\ref{9.4}), we find
  $$ 2a_{i+1}(n)+(n-1)a_i(n)=2a_{i+1}(n-2)+$$ $$2(2n-3)a_i(n-2)+2(n-1)(n-2)a_{i-1}(n-2)+$$
 \begin{equation}\label{9.6}
  \frac{1}{2}Coef[x^{m-i}] ((4x+3n-4)(x+\frac{n-4}{2})(x+\frac{n-6}{2})\cdot...\cdot(x+1)(x)).
   \end{equation}
   \newpage
Note that, polynomial $(4x+3n-4)(x+\frac{n-4}{2})(x+\frac{n-6}{2})\cdot...\cdot(x+1)x$ has degree $m+1.$ Therefore, in order to calculate $Coef[x^{m-i}]$ in (\ref{9.6}), we should choose, by all possible ways, in $m-i$ brackets (from $m+1$ ones) $x's,$ and in other $i+1$ brackets we choose linear forms of $n.$ Thus $\frac{1}{2}Coef[x^{m-i}]$ in (\ref{9.6}) is a polynomial $r_{i}(n)$ of degree $i+1.$ Further we use induction over $i$ with the formulas (\ref{9.2})-(\ref{9.3}) as the inductive base. Write (\ref{9.6}) in the form
  $$ 2(a_{i+1}(n)-a_{i+1}(n-2))=$$ $$2(2n-3)a_i(n-2)-(n-1)a_i(n)+$$
 \begin{equation}\label{9.7}
  2(n-1)(n-2)a_{i-1}(n-2)+r_i(n).
   \end{equation}
 By the inductive supposition, $a_{i-1}(n),$ $a_i(n)$ are polynomials of degree $2i-1$ and $2i+1$ respectively. Thus  $a_{i+1}(n)-a_{i+1}(n-2)$ is a polynomial of degree $2i+2.$ This means that $a_{i+1}$ is a polynomial of degree $2i+3.$\newline
 2) Let $n$ be odd. By (\ref{6.3}), for integer $x$ and $m=\frac{n-1}{2},$ we have
 $$ (2x+n-2)(a_0(n)x^m+a_1(n)x^{m-1}+...)=$$ $$2(x+n-1)(x+n-2)(a_0(n-2)x^{m-1}+a_1(n-2)x^{m-2}+...)+$$
\begin{equation}\label{9.8}
 (\frac{n-1}{2})!(4x+3n-4)\binom{x+\frac{n-3}{2}}{\frac{n-1}{2}}.
 \end{equation}
 Hence, comparing the coefficient of $x^{m+1}$ in both hand sides, we find
 $$a_0(n)=a_0(n-2)+2,\enskip n\geq3,\enskip a_0(1)=1.$$
 Thus $a_0(3)=3, a_0(5)=5,..., a_0(n)=n.$  \newline
 Furthermore, comparing the coefficient of $x^{m}$ in both hand sides in (\ref{9.8}), using the same arguments as in 1), we have
 $$a_1(n)=a_1(n-2)+\frac {7n^2-22n+19}{4},\enskip n\geq3,\enskip a_1(1)=0.$$
 Since $a_1(n)=\sum_{i=3,5,...,n}(a_1(i)-a_1(i-2)),$ then we find
 $$a_1(n)=\frac{1}{4}\sum_{i=3,5,...,n}(7i^2-22i+19)=\frac{1}{24}(7n^3-12n^2+5n).$$
  Finally, comparing the coefficient of $x^{m-i}$ in both hand sides of (\ref{9.8}), we find
$$ 2(a_{i+1}(n)-a_{i+1}(n-2))=$$ $$2(2n-3)a_i(n-2)-(n-2)a_i(n)+$$
\begin{equation}\label{9.9}
2(n-1)(n-2)a_{i-1}(n-2)+s_{i}(n),
\end{equation}
where
$$s_{i}(n)= Coef[x^{m-i}] ((4x+3n-4)(x+\frac{n-3}{2})(x+\frac{n-5}{2})\cdot...\cdot(x+1)x)$$
and, as in 1), the statement is proved by induction over $i. \blacksquare$
\newpage
  A few such polynomials are the following:\newline
 For odd $n:$
$$ U_0(n)=n,$$
$$ U_1(n)=\frac{1}{24}(n-1)n(7n-5),$$
$$U_2(n)=\frac{1}{640}(n-3)(n-1)n(29n^2-44n+7),$$
$$U_3(n)=\frac{1}{322560}(n-5)(n-3)(n-1)n(1581n^3-3775n^2+1587n+223); $$\newline
For even $n:$
$$ V_0(n)=\frac{1}{2}n,$$
$$ V_1(n)=\frac{1}{48}(n-2)n(7n-4),$$
$$V_2(n)=\frac{1}{3840}(n-4)(n-2)n(87n^2-98n+16),$$
$$V_3(n)=\frac{1}{645120}(n-6)(n-4)(n-2)n(1581n^3-2686n^2+936n+64). $$

 \begin{proposition}\label{P2}
  \begin{equation}\label{9.10}
a_i(n)\equiv \begin{cases}
r_{i}(n), \;\;if\;\;n\;\;is\;\;even,\\\\s_{i}(n),\;\;if\;\; n\;\; is\;\; odd.\end{cases}\pmod2
\end{equation}
 \end{proposition}
  \bfseries Proof. \mdseries The proposition follows from ({9.7}), ({9.9}) and Theorem 2. $\blacksquare$ \newline \indent Finally, note that, from ({8.7})-({8.8}) follow the following homogeneous recursions for the coefficients of $P_n(x).$
\begin{theorem}\label{t7} For odd $n\geq3$ and $i\geq0,$
 \begin{equation}\label{9.11}
(m-i)a_i(n)=na_i(n-1)+\sum_{j=0}^{i-1}(-1)^{i-j+1}\binom{m-j}{m-i-1}a_j(n).
\end{equation}
For even $n\geq4$ and $i\geq0,$
$$(n-2i-1)a_i(n)=\frac{n}{2}a_i(n-1)+$$
\begin{equation}\label{9.12}
2\sum_{j=0}^{i-1}(-1)^{i-j+1}(m\binom{m-j}{m-i}-\binom{m-j}{m-i-1})a_j(n).
\end{equation}
\end{theorem}
\newpage
\section{Arithmetic proof of the integrality $P_n(x)$ in integer points}
From Theorem 2 we conclude that the polynomial $P_n(x)$ takes integers values for integer $x=k.$ Here we give an independent arithmetic proof of this fact, using the explicit expression (\ref{1.6}).
It is well known (cf. \cite{5}, Section 8, Problem 87) that, if a polynomial $P(x)$ of degree $m$ takes integer values for $x=0,1,...,m,$ then it takes integer values for every integer $x.$ Since, as we proved at the end of Section 6, $\deg {P_n(k)}=\lfloor\frac{n-1}{2}\rfloor,$ then we suppose that $0\leq k\leq \lfloor\frac{n-1}{2}\rfloor.$ Moreover, from the results of Section 3, $P_n(0)$ and $P_n(1)$ are integers (in the case when $n+1$ is an odd prime, $P_n(1)=(2^n-1)(\frac{n}{2})!/(n+1)$ is integer, since $2^n-1\equiv0 \pmod{n+1},$ while in the case when $n+1$ is an odd composite number, no divisor exceeds $\frac {n+1}{3},$ therefore, $(\frac{n}{2})!\equiv0 \pmod{n+1}).$ Thus we can suppose that
\begin{equation}\label{10.1}
2\leq k\leq \lfloor\frac{n-1}{2}\rfloor.
\end{equation}
Suppose that $n$ is \slshape even \upshape  (the case of odd $n$ is considered quite analogously).
Let $p$ be a prime. Denote the maximal power of $p$ dividing $n$ by $[n]_p.$ We say that, for integer $l,\enskip h,$ the fraction $\frac{l}{h}$ is \slshape $p$-integer,\upshape \enskip if $[l]_p-[h]_p\geq0.$\newline
A) Firstly, we show that, for $n\geq4,$ $P_n(k)$ is $2$-integer. Indeed, $2k+n-1$ is odd, while $4k+3n-4$ is even. Therefore, by (\ref{6.4}), using a trivial induction, we see that $P_{n}(k)$ is $2$-integer. \newline
\indent Further we use the explicit formula (\ref{1.6}) of Theorem \ref{t1}. \newline
B) Let $p$ be an odd prime divisor of $\binom{n+2k-2}{k-1}$ which does not coincide with any factor of the product $(n+2k-1)(n+2k-2)...(n+k).$ Thus $p$ could divide one or several \slshape composite \upshape \enskip factors of this product. Therefore, the following condition holds
\begin{equation}\label{10.2}
 3\leq p\leq \frac{n+2k-1}{3}.
 \end{equation}
 Let us show that
 $$a(n;k):=\frac {\binom{\frac{n}{2}+k-1}{k}}{\binom{n+2k-1}{k}}(\frac{n-2}{2})!=$$
\begin{equation}\label{10.3}
2^{-k}\frac{(n+2k-2)(n+2k-4)\cdot...\cdot n}{(n+2k-1)(n+2k-2)...(n+k)}(\frac{n-2}{2})!
 \end{equation}
is $p$-integer and, consequently, $P_n(k)$ is $p$-integer. \newline
Let $k\geq3$ be even. Then, after a simplification, we have
$$2^ka(n;k)=\frac{(n+k-2)(n+k-4)\cdot...\cdot n}{(n+2k-1)(n+2k-3)\cdot...\cdot(n+k+1)}(\frac{n-2}{2})!,$$
or
\begin{equation}\label{10.4}
2^{\frac{k}{2}}a(n;k)=\frac{(\frac{n+k-2}{2})!}{(n+2k-1)(n+2k-3)...(n+k+1)}
 \end{equation}
 \newpage
 We distinguish several cases.\newline
\indent Case a) For $t\geq2,$ let $p^t$ divide at least one factor of the denominator. Then  $p\leq (n+2k-1)^{\frac{1}{t}}$ Let us show that $p\leq  \frac{n+k-2}{2t}.$ We should show that $n+2k-1\leq (\frac{n+k-2}{2t})^t,$ or, since, by ({10.1}), $k\leq\frac{n-2}{2},$ it is sufficient to show that $\frac{3}{2}(n+k-2)\leq (\frac{n+k-2}{2t})^t,$ or $(2t)^{\frac{t}{t-1}}\leq (\frac{2}{3})^{\frac{1}{t-1}}(n+k-2).$ Since $(\frac{2}{3})^{\frac{1}{t-1}}\geq\frac{2}{3},$
it is sufficient to prove that $(2t)^{\frac{t}{t-1}}\leq \frac{2}{3}(n+k-2).$ Note that $e^t<p^t\leq n+2k-2,\enskip t\leq\ln(n+2k-2).$ Therefore we find $(2t)^{\frac{t}{t-1}}\leq (2\ln (n+2k-2))^2.$ Furthermore, note that, if $n\geq152,$ then $\ln^2n<\frac{n}{6}.$ Thus $(2t)^{\frac{t}{t-1}}\leq \frac{2}{3}(n+k-2).$ It is left to add that up to $n=161$ we verified that the polynomials $P_n(k)$ have integer coefficients and, consequently, is integer-valued.\newline
\indent Case b) Let $p$ divide only one factor of the denominator. Then, in view of ({10.1}) and ({10.2}), $p\leq\frac{n+2k-1}{3}\leq \frac{n+k-2}{2}$ and, by ({10.4}), $a(n;k)$ is $p$-integer.\newline
\indent Case c) Let $p$ divide exactly $l$ factors of the denominator. Then $p\leq \frac{(n+2k-1)-(n+k+1)}{l}=\frac{k-2}{l},$ and, since, by ({10.1}), $n\geq2k+2,$ we conclude that $\frac{n+k-2}{2}\geq\frac {3k}{2}\geq k-2\geq lp.$ Hence, by ({10.4}), $a(n;k)$ is $p$-integer.\newline
It is left to notice that the case of odd $k$ is considered quite analogously. \newline
C) Suppose that, as in B), $k\geq2$ is even. Let $p$ be an odd prime divisor of $\binom{n+2k-1}{k}$ which coincides with some factor of the product $(n+2k-1)(n+2k-3)...(n+k+1).$ In this case the fraction ({10.4}) is not integer. Thus in order to prove that $P_n(k)$ is $p$-integer, we should prove that $T_n(k)$ (\ref{5.1}) is $p$-integer. By the condition, $p$ has form
\begin{equation}\label{10.5}
p=n+2k-1-2r,\enskip 0\leq r\leq \frac {k-2}{2}.
 \end{equation}
According to (\ref{5.1}) and ({10.5}), we should prove that
$$\sum_{j=0}^{n-1} 2^{j}\binom{n+2k-j-2}{k-1}=$$
\begin{equation}\label{10.6}
\sum_{j=0}^{n-1} 2^{j}\binom{p+2r-1-j}{k-1}\equiv0\pmod{p},
 \end{equation}
or
$$A(n,r,k):=$$ $$\sum_{j=0}^{n-1} 2^{j}(j-(2r-1))(j+1-(2r-1))...(j+k-2-(2r-1))\equiv0\pmod{p}.$$
Note that, since $n-2r=p-2k+1,$ then we have
$$\sum_{j=0}^{n-1}x^{j+k-2-(2r-1)}=$$ $$(x^{n+k-2r-1}-x^{k-2r-1})(x-1)^{-1}=(x^{p-k}-x^{k-2r-1})(x-1)^{-1}.$$
\newpage
Therefore,
$$A(n,r,k)=2^{2r}\sum_{j=0}^{n-1}(x^{j+k-2-(2r-1)})^{(k-1)}\mid_{x=2}=$$ $$2^{2r}((x^{p-k}-x^{k-2r-1})(x-1)^{-1})^{(k-1)}\mid_{x=2}.$$
Thus we should prove that
\begin{equation}\label{10.7}
((x^{p-k}-x^{k-2r-1})(x-1)^{-1})^{(k-1)}\mid_{x=2}\equiv0\pmod{p},
 \end{equation}
or, using the Leibnitz formula,
$$\sum_{j=0}^{k-1}(-1)^{k-j-1}\binom{k-1}{j}(k-j-1)!(p-k)(p-k-1)...(p-k-j+1)2^{p-k-j}\equiv $$ $$\sum_{j=0}^{k-1}(-1)^{k-j-1}\binom{k-1}{j}(k-j-1)!(k-2r-1)(k-2r-2)...(k-2r-j)2^{k-2r-j-1}\pmod p.$$
Since $2^{p-1}\equiv1\pmod p,$ then we should prove the identity
$$\sum_{j=0}^{k-1}(-1)^{k-j-1}\binom{k-1}{k-j-1}(k-j-1)!(p-k)(p-k-1)...(p-k-j+1)\mid_{p=0}2^{-k-j+1}= $$ $$\sum_{j=0}^{k-1}(-1)^{k-j-1}\binom{k-1}{k-j-1}(k-j-1)!(k-2r-1)(k-2r-2)...(k-2r-j)2^{k-2r-j-1},$$
or, after simple transformations, the identity
\begin{equation}\label{10.8}
\sum_{j=0}^{k-1}\binom{k+j-1}{j}2^{-j}=2^{2k-2r-2}\sum_{j=0}^{k-1}(-1)^{j}\binom{k-2r-1}{j}2^{-j}.
 \end{equation}
It is known ((\cite{7},\enskip Ch.1, problem 7), that
 $$\sum_{i=0}^{n}\binom{2n-i}{n}2^{i-n}=2^n.$$
 Putting $n-i=j,$ we have
 $$\sum_{j=0}^{n}\binom{n+j}{n}2^{-j}=\sum_{j=0}^{n}\binom{n+j}{j}2^{-j}=2^n.$$
 Therefore, the left hand side in ({10.8}) is $2^{k-1}$ and it is left to prove that
 $$\sum_{j=0}^{k-1}(-1)^{j}\binom{k-2r-1}{j}2^{k-j}=2^{2r+1}.$$
 We have
 $$\sum_{j=0}^{k-1}(-1)^{j}\binom{k-2r-1}{j}2^{k-j}=\sum_{j=0}^{k-2r-1}(-1)^{j}\binom{k-2r-1}{j}2^{k-j}=$$ $$2^{2r+1}\sum_{j=0}^{k-2r-1}(-1)^{j}\binom{k-2r-1}{j}2^{k-2r-1-j}=2^{2r+1}(2-1)^{k-2r-1}=2^{2r+1}$$
 \newpage
 and we are done. The case of odd $k\geq3$ is considered quite analogously. So, formulas ({10.4})-({10.5}) take the form
 $$2^{\frac{k-1}{2}}a(n;k)=\frac{(\frac{n+k-1}{2})!}{(n+2k-1)(n+2k-3)\cdot...\cdot(n+k)},$$
\begin{equation}\label{10.9}
 p=n+2k-2r-1,\enskip 0\leq r\leq \frac {k-1}{2},
 \end{equation}
 and, for odd $k,$ the proof reduces to the same congruence ({10.6}).
  $\blacksquare$

\section{Representation of $P_n(x)$ in basis $\{\binom {x}{i}\}$}
The structure of explicit formula (\ref{1.6}) allows us to conjecture that the coefficients of $P_n(x)$ in basis $\{\binom {x}{i}\}$ possess simpler properties.
A process of expansion of a polynomial $P(x)$ in the binomial basis is indicated in \cite{5} in a solution of Problem 85: ``Functions $1,\enskip x,\enskip x^2,...,x^n$ one can consecutively express in the form of linear combinations with the constant coefficients of $ 1,\enskip \frac {x}{1},\enskip \frac{x(x-1)}{2},...,\frac{x(x-1)...(x-n+1)}{n!}."$ Therefore,
$$P(x)=b_0\binom {x}{m}+b_1\binom{x}{m-1}+...+b_{m-1}\binom{x}{1}+b_m,$$
where $b_0,\enskip b_1,...,b_m$ are defined from the equations
$$P(0)=b_m,$$
$$P(1)=b_m+\binom{1}{1}b_{m-1},$$
$$P(2)=b_m+\binom{2}{1}b_{m-1}+\binom{2}{2}b_{m-2},$$
$$. \enskip. \enskip.\enskip.\enskip.\enskip.\enskip.\enskip.\enskip.\enskip.\enskip.\enskip.\enskip.\enskip.
\enskip.\enskip.\enskip.\enskip.\enskip.\enskip.\enskip.$$
\begin{equation}\label{11.1}
P(m)=b_m+\binom{m}{1}b_{m-1}+...+\binom{m}{m}b_0.
 \end{equation}
This process one can simplify in the following way. In the identity
$$n^x=(1+(n-1))^x=$$ $$1+(n-1)\binom{x}{1}+(n-1)^2\binom{x}{2}+...+(n-1)^x\binom{x}{x}=$$
$$n^0+(n-n^0)\binom{x}{1}+(n-n^0)^2\binom{x}{2}+...+(n-n^0)^x\binom{x}{x}  $$
we can evidently replace powers $n^j,\enskip j=0,...,x,$ by the arbitrary numbers $a_j,\enskip j=0,...,x.$ Thus we have a general identity
$$a_x=a_0+(a_1-a_0)\binom{x}{1}+$$ $$(a_2-2a_1+a_0)\binom{x}{2}+
(a_3-3a_2+3a_1-a_0)\binom{x}{3}+...+$$
\newpage
\begin{equation}\label{11.2}
(a_x-\binom{x}{1}a_{x-1}+
\binom{x}{2}a_{x-2}-...+(-1)^x\binom{x}{x}a_0)\binom{x}{x}.
 \end{equation}
 Essentially, we quickly obtained a special case of the so-called ``Newton's forward difference formula" (cf. \cite {10}).
  Here, put $a_j=P(j),\enskip j=0,...,m,$ and, firstly, consider values $0\leq x\leq m.$ Since $\binom {x}{l}=0$ for $l>m,$ then we obtain the required representation under the condition $0\leq x\leq m:$
 $$P(x)=P(0)+(P(1)-P(0))\binom{x}{1}+$$$$(P(2)-2P(1)+P(0))\binom{x}{2}+...+
(P(m)-\binom{m}{1}P(m-1)+$$
\begin{equation}\label{11.3}
\binom{m}{2}P(m-2)-...+(-1)^m\binom{m}{m}P(0))\binom{x}{m}.
 \end{equation}
It is left to note that, since a polynomial of degree $m$ is fully defined by its values in $m+1$ points $0,1,...,m,$ then (\ref{11.3}) is the required representation for all $x.$\newline
 \indent So, for the considered polynomials $\{P_n(x)\},$ we have
$$P_1=1,$$
$$P_2=1,$$
$$P_3=3\binom {x}{1}+4,$$
$$P_4=2\binom {x}{1}+4,$$
$$P_5=10\binom {x}{2}+30\binom {x}{1}+32,$$
$$P_6= 6\binom {x}{2}+22\binom {x}{1}+32,$$
$$P_7=42\binom {x}{3}+196\binom {x}{2}+378\binom {x}{1}+384,$$
$$P_8=24\binom {x}{3}+128\binom{x}{2}+296\binom{x}{1}+384,$$
$$P_9=216\binom{x}{4}+1368\binom {x}{3}+3816\binom {x}{2}+6120\binom {x}{1}+6144,$$
$$P_{10}=120\binom{x}{4}+840\binom {x}{3}+2664\binom {x}{2}+5016\binom {x}{1}+6144,$$
$$P_{11}=1320\binom{x}{5}+10560\binom {x}{4}+38544\binom {x}{3}+84480\binom {x}{2}+122760\binom{x}{1}+122880,$$
$$P_{12}= 760\binom{x}{5}+ 6240\binom {x}{4}+25152\binom {x}{3}+62112\binom {x}{2}+103920\binom{x}{1}+122880.$$
\newpage
 \section{On coefficients of $P_n(x)$ in basis $\{\binom {x}{i}\}$}
 Let
 \begin{equation}\label{12.1}
 P_n(x)=b_0(n)\binom {x}{m}+b_1(n)\binom{x}{m-1}+...+b_{m-1}(n)\binom{x}{1}+b_m(n),
  \end{equation}
where $m=\lfloor\frac{n-1}{2}\rfloor.$\newline
 Since, for integer $k,$ we have the explicit formula for $P_n(k)$ (\ref{1.6}), then, according to (\ref{11.3}), we have the following explicit formula for $b_i(n),\enskip i=0,...,m:$
  \begin{equation}\label{12.2}
 b_i(n)=\sum_{k=0}^{m-i}(-1)^{m-i-k}\binom{m-i}{k}P_n(k).
  \end{equation}
  Let
  $$P_n(x)=\sum_{j=0}^m a_j(n){x}^{m-j}.$$
  Then
  \begin{equation}\label{12.3}
  b_i(n)=\sum_{j=0}^{m}a_j(n)\sum_{k=0}^{m-i}(-1)^{m-i-k}k^{m-j}\binom{m-i}{k}.
  \end{equation}
  Since the $l$-th difference of $f(x)$ is (cf. \cite{1}, formula 25.1.1)
   $$\Delta^l f(x)=\sum_{k=0}^{l}(-1)^{l-k}\binom{l}{k}f(x+k),$$
   then one can write (\ref{12.3})in the form
  $$b_i(n)=\sum_{j=0}^{m}a_j(n)\Delta^{m-i}x^{m-j}\mid_{x=0}.$$
Here the summands corresponding to $ j>i,$ evidently, equal 0. Therefore, we have
\begin{equation}\label{12.4}
  b_i(n)=\sum_{j=0}^{i}a_j(n)\Delta^{m-i}x^{m-j}\mid_{x=0}.
\end{equation}
\begin{theorem}\label{t8}
 For $n\geq1,$ we have
\begin{equation}\label{12.5}
 b_0(n)=\begin{cases}
n(\frac{n-1}{2})!, \;\;if\;\;n\;\;is\;\;odd,\\\\(\frac{n}{2})!,\;\;if\;\; n\;\; is\;\; even;\end{cases}
\end{equation}

\begin{equation}\label{12.6}
b_1(n)=\begin{cases}
\frac{1}{6}n(5n-7)(\frac{n-1}{2})!, \;\;if\;\;n\;\;is\;\;odd,\\\\\frac{1}{6}(5n-8)(\frac{n}{2})!,\;\;if\;\; n\;\; is\;\; even.\end{cases}
\end{equation}
In general, for a fixed $i,$ $b_i(n)=(m-i)!Y_i(n),$ if $n$ is odd, and $b_i(n)=(m-i)!Z_i(n),$ if $n$ is even, where $Y_i,\enskip Z_i$ are polynomials in $n$ of degree $2i+1.$
\end{theorem}
\newpage
 \bfseries Proof. \mdseries Note that the Stirling number of the second kind $S(n,m)$ is connected with the $m$-th difference of $\Delta^m x^n\mid_{x=0}$ in the following way (see \cite{1}, formulas 24.1.4)
\begin{equation}\label{12.7}
S(n,m)m!=\sum_{k=0}^{m}(-1)^{m-k}\binom{m}{k}k^n=\Delta^m x^n\mid_{x=0}.
  \end{equation}
  In particular, since $S(m,m)=1,$ $S(m+1,m)=\binom{m+1}{2},$ then
  $$\Delta^{m}x^{m}\mid_{x=0}=m!$$
  and
\begin{equation}\label{12.8}
     \Delta^{m}x^{m+1}\mid_{x=0}=\frac{m}{2}(m+1)!
\end{equation}
Therefore, by (\ref{12.4}),
  $$b_0(n)=m!a_0(n),$$\newline
  $$b_1(n)=\frac{m-1}{2}m!a_0(n)+(m-1)!a_1(n),$$\newline
  and, by (\ref{9.2})-(\ref{9.3}) (where $m=\lfloor\frac{n-1}{2}\rfloor$), we find formulas  (\ref{12.5})-(\ref{12.6}).$\blacksquare$
  Further we need lemma.
  \begin{lemma}\label{L6}
  $S(n+k,n)$ is a polynomial in $n$ of degree $2k.$
  \end{lemma}
\bfseries Proof. \mdseries For $k\geq1,$ denote
\begin{equation}\label{12.9}
Q_k(n)=S(n+k,n).
\end{equation}
Note that, since $S(n,n)=1,$ then $Q_0(n)=1.$ Further, since $S(n,0)=\delta_{n,0},$ then, for $k\geq1,$ $Q_k(0)=0.$
From the main recursion for $S(n,m)$ which is $S(n,m)=mS(n-1,m)+S(n-1,m-1),$ we have
\begin{equation}\label{12.10}
Q_k(n)-Q_k(n-1)=nQ_{k-1}(n).
\end{equation}
and, in view of $Q_k(0)=0,$ we find the recursion
\begin{equation}\label{12.11}
Q_0(n)=1,\enskip Q_k(n)=\sum_{i=1}^niQ_{k-1}(i).
\end{equation}
Using a simple induction, from (\ref{12.11}) we obtain the lemma. $\blacksquare$
\begin{remark} The list of polynomials $\{Q_k(n)\}$
$$Q_0=1,$$
$$Q_1=\frac{1}{2}{n(n+1)}, $$
$$Q_2=\frac{1}{24}{n(n+1)(n+2)(3n+1)},$$
$$ Q_3=\frac{1}{48}{n^2(n+1)^2(n+2)(n+3)},$$
\newpage
$$ Q_4=\frac{1}{5760}{n(n+1)(n+2)(n+3)(n+4)(15n^3+30n^2+5n-2)}, etc.$$\newline
It could be proven that the sequence of denominators coincides with $A053657$ \cite{8}, such that the denominator of $Q_k(n)$ is
$\prod p^{\sum_{j\geq0}\lfloor \frac{k}{(p-1)p^j}\rfloor},$ where the product is over all primes.
\end{remark}
\indent Note that from (\ref{12.4}) and (\ref{12.7}) we find
\begin{equation}\label{12.12}
 b_i(n)=(m-i)!\sum_{j=0}^{i}a_j(n)S(m-j,m-i), \enskip m=\lfloor\frac{n-1}{2}\rfloor.
  \end{equation}
 Since, by Lemma \ref{L5}, $S(m-j,m-i)$ is a polynomial in $n$ of degree $2((m-j)-(m-i))=2(i-j),$ while, by Theorem \ref{t6}, $a_j(n)$ is a polynomial of degree $2j+1,$ then $a_j(n)S(m-j,m-i)$ is a polynomial of degree $2i+1.$ Thus
 $\sum_{j=0}^{i}a_j(n)S(m-j,m-i)$ is a polynomial of degree $2i+1.$ This completes the proof. $\blacksquare$ \newline
 The first polynomials $Y_i(n),\enskip Z_i(n)$ are
 $$Y_0=n,$$
 $$ Y_1=\frac{1}{12}(n-1)n(5n-7),$$
 $$Y_2=\frac{1}{480}(n-3)(n-1)n(43n^2-168n+149),$$
 $$Y_3=\frac{1}{13440}(n-5)(n-3)(n-1)n(177n^3-1319n^2+3063n-2161);$$
 $$Z_0=\frac{n}{2},$$
 $$ Z_1=\frac{1}{24}(n-2)n(5n-8),$$
 $$Z_2=\frac{1}{960}(n-4)(n-2)n(43n^2-182n+184),$$
 $$Z_3=\frac{1}{26880}(n-6)(n-4)(n-2)n(3n-8)(59n^2-306n+352)).$$
Finally, we prove the following attractive result.
\begin{theorem}\label{t9}
 $1)$ For odd $n,$ $b_j(n)/n, j=0,...,m-1,$ are integer. Moreover, for $n\geq3,$
 \begin{equation}\label{12.13}
 b_i(n)=n(b_{i}(n-1)+b_{i-1}(n-1)),\enskip i=1,...,m-1.
  \end{equation}
 $2)$ For even $n\geq4,$
 \begin{equation}\label{12.14}
 2b_i(n)=b_{i}(n-1)+b_{i-1}(n-1)+m!\binom{m}{i},\enskip i=1,...,m-1.
  \end{equation}
 \end{theorem}
 \newpage
 \bfseries Proof. \mdseries 1) According to (\ref{12.2}), we should prove that for odd $n\geq3,$
$$\sum_{k=0}^{m-i}(-1)^{m-i-k}\binom{m-i}{k}P_n(k)=$$
 $$n(\sum_{k=0}^{m-i}(-1)^{m-i-k}\binom{m-i}{k}P_{n-1}(k)+$$ $$\sum_{k=0}^{m-i-1}(-1)^{m-i-k-1}\binom{m-i-1}{k}P_{n-1}(k)),\enskip i=1,2,...,m-1,$$
 or, putting $m-i=t,$
 $$\sum_{k=0}^{t}(-1)^{k} \binom{t}{k}P_n(k)=n(\sum_{k=0}^{t}(-1)^{k}\binom{t}{k}P_{n-1}(k)-$$ $$\sum_{k=0}^{t}(-1)^{k}\binom{t-1}{k}P_{n-1}(k)),\enskip t=1,2,...,m-1,$$
 or, finally, for $t=1,...,\frac{n-3}{2},$
 \begin{equation}\label{12.14}
\sum_{k=1}^{t}(-1)^{k-1} (\binom{t}{k}P_n(k)-\binom {t-1}{k-1}nP_{n-1}(k))=P_n(0).
 \end{equation}
To prove ({12.15}), note that , by (\ref{8.7}), $nP_{n-1}(k)=P_n(k)-P_n(k-1).$ Hence,
$$ \binom{t}{k}P_n(k)-\binom {t-1}{k-1}nP_{n-1}(k)=$$ $$P_n(k)(\binom{t}{k}-\binom {t-1}{k-1})+\binom {t-1}{k-1}P_n(k-1)=$$ $$\binom{t-1}{k}P_n(k)+\binom {t-1}{k-1}P_n(k-1). $$
Thus the summands of ({12.15}) are
$$(-1)^{k-1} (\binom{t}{k}P_n(k)-\binom {t-1}{k-1}nP_{n-1}(k))=$$ $$(-1)^{k-1} \binom{t-1}{k}P_n(k)-(-1)^{k-2}\binom {t-1}{k-1}P_{n}(k-1),$$
and the summing gives
$$\sum_{k=1}^{t}(-1)^{k-1} (\binom{t}{k}P_n(k)-\binom {t-1}{k-1}nP_{n-1}(k))=$$ $$(-1)^{k-1} \binom{t-1}{k}P_n(k)\mid_{k=t}-(-1)^{k-2}\binom {t-1}{k-1}P_{n}(k-1)\mid_{k=1}=P_n(0). $$
2) Analogously, the proof of ({12.14}) reduces to proof of the following equality for $t=1,2,...,m-1:$
\begin{equation}\label{12.16}
\sum_{k=0}^{t+1}(-1)^{k} (2\binom{t}{k}P_n(k)+\binom {t}{k-1}P_{n-1}(k))=(-1)^tm!\binom{m}{t}.
 \end{equation}
 \newpage
Note that, by (\ref{8.11}),
$$(-1)^{k} (2\binom{t}{k}P_n(k)+\binom {t}{k-1}P_{n-1}(k))=$$
\begin{equation}\label{12.17}
 (-1)^k\binom{t}{k}P_{n-1}(k+1)-$$ $$(-1)^{k-1}\binom {t}{k-1}P_{n-1}(k)+(-1)^k\binom{t}{k}\binom{k+m}{m}m!
 \end{equation}
 Since
 $$\sum_{k=0}^{t+1}((-1)^k\binom{t}{k}P_{n-1}(k+1)-(-1)^{k-1}\binom {t}{k-1}P_{n-1}(k))=$$
 $$(-1)^{k} \binom{t}{k}P_{n-1}(k)\mid_{k=t+1}-(-1)^{k-1}\binom {t}{k-1}P_{n-1}(k)\mid_{k=0}=0,$$
 then, by (\ref{12.16})-(\ref{12.17}), the proof reduces to the known combinatorial identity
 $$\sum_{k=0}^{t}(-1)^k\binom{t}{k}\binom{k+m}{m}=(-1)^t\binom{m}{t}, \enskip t=1,...,m-1$$
 (see \cite{7}, Ch.1, formula (8) with $p=0$ up to the notations). $\blacksquare$ \newline

\indent In conclusion, note that we have an interesting voyage from the, until now, unproved conjecture 7), when the primes in 7) are replaced by the odd numbers in Theorem \ref{t9}.\newline

{

\end {document}